\documentclass{amsart}
\usepackage{amssymb,amsmath}

\newtheorem{Theorem}{Theorem}[section]

\newtheorem{Proposition}[Theorem]{Proposition}
\newtheorem{Corollary}[Theorem]{Corollary}

\theoremstyle{definition}
\newtheorem{Remark}[Theorem]{Remark}

\newtheorem{Example}[Theorem]{Example}
\newtheorem{Conjecture}[Theorem]{Conjecture}

\numberwithin{equation}{section}

\begin{document}
\title{A Common Formula for Certain Generalized Hankel Transforms}
\author{Mario Garcia-Armas}
\address{Department of Mathematics, University of British Columbia, Vancouver, \newline \indent %
BC V6T 1Z2, Canada}
\email{marioga@math.ubc.ca}

\begin{abstract}
In this paper, we study the generalized Hankel transform of the family of sequences satisfying the recurrence relation $a_{n+1} = \bigl(\alpha + \frac{\beta}{n+\gamma}\bigr) a_n$. We find a connection between some well known formulas that had previously arisen in literature in dissimilar settings. 
\end{abstract}

\maketitle

\section{Introduction}
\label{SecnIntro}

We recall some terminology from the theory of Hankel matrices. Given a sequence $(a_n)_{n=0}^{\infty}$, we consider the doubly-indexed sequence of Hankel matrices $H_n^{(k)}$, $n=1,2,\ldots$, $k=0,1,\ldots$, defined by
\begin{equation} \label{Matrix_H_n_k}
H_n^{(k)} = \left( \begin{array} {ccccc} a_k&a_{k+1}&a_{k+2}&\ldots & a_{k+n-1}\\
a_{k+1}&a_{k+2}&a_{k+3}&\ldots&a_{k+n}\\
a_{k+2}&a_{k+3}&a_{k+4}&\ldots&a_{k+n+1}\\
\vdots&\vdots&\vdots&\ddots&\vdots\\
a_{k+n-1}&a_{k+n}&a_{k+n+1}&\ldots&a_{k+2n-2}
\end{array} \right).
\end{equation}

The \emph{(generalized) Hankel transform} of $(a_n)_{n=0}^{\infty}$ is the doubly-indexed sequence of determinants $d_n^{(k)} = \det H_n^{(k)}$ (for a similar treatment, see Garcia Armas and Sethuraman \cite{Garcia_Seth} and Tamm \cite{Tamm}). It is important to mention that several authors refer to the Hankel transform only as the sequence $d_n = \det H_n^{(0)}$ (see, for example, Chamberland and French \cite{Chamb_French} and Layman \cite{Lay}).

The theory of Hankel matrices have beautiful connections with many areas of mathematics, physics and computer science (see, for example, Desainte-Catherine and Viennot \cite{DCV}, Garcia Armas and Sethuraman \cite{Garcia_Seth}, Tamm \cite{Tamm}, Vein and Dale \cite{Vein_Dale}). Although Hankel determinants had been previously studied (see, for example, Aigner \cite{Aigner}), the term \emph{Hankel transform} was introduced and first studied by J. W. Layman in \cite{Lay}. Several later studies of Hankel transforms of sequences have appeared in literature (see, for example, Chamberland and French \cite{Chamb_French}, Cvetkovi\' c et al. \cite{CRI}, Egecioglu et al. \cite{ERR_Riemann, ERR}, French \cite{French}, Spivey and Steil \cite{SpSt}).

In the evaluation of Hankel determinants, several techniques have proved to be useful. For an extended set of tools, as well as a significant bibliography, please refer to Krattenthaler \cite[Sec. 2.7]{Kratt1} (also \cite[Sec. 5.4]{Kratt2}) and Vein and Dale \cite{Vein_Dale}. 


In this note, we study the generalized Hankel transform of a sequence $(a_n)_{n=0}^{\infty}$ satisfying
\begin{equation} \label{Rel_def_a_n}
a_{n+1} = \biggl(\alpha + \frac{\beta}{n+\gamma}\biggr) a_n, \qquad \forall n \geq 0
\end{equation}
for some complex numbers $\alpha, \beta, \gamma$. We would like to emphasize that we do not claiming originality for the evaluation of such Hankel transform; it can be achieved using, e.g., the work of Krattenthaler \cite{Kratt1}. Our main goal is to illustrate the hidden connection between several Hankel transform evaluations that have been previously studied in quite independent settings.

\section{Basic Computations}
\label{SecnBasicComp}

Throughout this paper, we consider a finite product $\prod_{s=a}^{b} c_s = 1$ when $b<a$.

Let $(a_n)_{n=0}^{\infty}$ be a sequence satisfying \eqref{Rel_def_a_n}. For simplicity, we may assume that $a_0 = 1$. It is quite obvious that the results we derive in this section can be straightforwardly extended to the general case.

Note that for every $m \geq n \geq 0$ we can write
\begin{equation} \label{a_n_product_form}
a_m = a_n\prod_{s=n+1}^m \biggl(\frac{\alpha (s+\gamma-1) + \beta}{s+\gamma-1} \biggr) = a_n\prod_{s=1}^{m-n} \biggl(\frac{\alpha (s+n+\gamma-1) + \beta}{s+n+\gamma-1} \biggr).
\end{equation}

Let us consider the Hankel matrix $H_n^{(k)} = (m_{ij})_{1\leq i,j \leq n}$, where $m_{ij} = a_{i+j+k-2}$. From Equation \eqref{a_n_product_form} we obtain the relation
\begin{equation}
m_{ij} = a_{j+k-1} \prod_{s=1}^{i-1} \biggl(\frac{\alpha (s+j+k+\gamma-2) + \beta}{s+j+k+\gamma-2} \biggr).
\end{equation}

\begin{Proposition}
\begin{equation} \label{Ugly_formula_d_n_k}
d_n^{(k)} = a_k a_{k+1} \ldots a_{k+n-1}\ \frac{\displaystyle\prod_{i=1}^{n-1} i! \: \bigl[\alpha(i-1) - \beta \bigr]^{n-i}}{\displaystyle\prod_{i=1}^{n-1} \bigl[i+k+\gamma-1\bigr]^i \: \bigl[2n-i+k+\gamma-2\bigr]^i}.
\end{equation}
\end{Proposition}
\begin{proof}
This formula could be derived, after some work, from Krattenthaler \cite[Thm. 26, Eq. 3.11]{Kratt1}.

Alternatively, a different (and more elementary) approach will be outlined. Let $M_n(a,b,c)$ denote the $n \times n$ matrix with entries
\begin{equation}
m'_{ij} = \prod_{s=1}^{i-1} \frac{a(s+j)+b}{s+j+c}, \qquad i,j = 1,2,\ldots,n;
\end{equation}
where $a, b, c$ are complex numbers. Then, it can be proved that the following recurrence relation holds:
\begin{equation} \label{Recurrence_dets_M_n}
\det M_n(a,b,c)  = \frac{(n-1)! \: (ac-b)^{n-1}}{\displaystyle \prod_{i=1}^{n-1} \bigl( i+1+c \bigr) \bigl( i+2+c \bigr)} \; \det M_{n-1}(a,a+b,c+2).
\end{equation}
This can be achieved after the following sequence of steps: 
\begin{itemize}
\item Subtracting the $(k-1)$-st column of $M_n(a,b,c)$ from the $k$-th one for $k = n, n-1, \ldots, 2$, in that order.
\item Expanding by the minors of the first row.
\item Taking common factors from rows and columns out of the determinant of the lower right $(n-1) \times (n-1)$ block (These common factors form the fraction on the left of the RHS of Equation \eqref{Recurrence_dets_M_n}).
\end{itemize}


The proposition follows easily from Equation \eqref{Recurrence_dets_M_n} after the substitution
\begin{equation*}
(a,b,c) := (\alpha, \alpha \gamma + \alpha k -2 \alpha +\beta, \gamma + k - 2).
\end{equation*}
\end{proof}

Let us find a nicer expression for $d_n^{(k)}$. We make the following observation.

\begin{Remark} \label{Remark_d_n_k_null}
It is easily seen that $d_n^{(k)} = 0$ implies $d_n^{(k+1)} = 0$. Indeed, the numerator of the fraction on the right of \eqref{Ugly_formula_d_n_k} does not depend on $k$ and on the other hand, if $a_i = 0$ for some $i \in \{k,\ldots,k+n-1\}$, then Equation \eqref{Rel_def_a_n} yields $a_j = 0$ for all $j \geq i$.
\end{Remark}

Suppose now that we have $d_n^{(j)} \neq 0$ for some $j \geq 0$. Using Equations \eqref{a_n_product_form} and \eqref{Ugly_formula_d_n_k} together, we compute the ratio
\begin{eqnarray}
\frac{d_n^{(j+1)}}{d_n^{(j)}} &=& \frac{a_{j+n}}{a_j}\times \frac{\displaystyle\prod_{i=1}^{n-1} \bigl[i+j+\gamma-1\bigr]^i \: \bigl[2n-i+j+\gamma-2\bigr]^i}{\displaystyle\prod_{i=1}^{n-1} \bigl[i+j+\gamma\bigr]^i \: \bigl[2n-i+j+\gamma-1\bigr]^i} \nonumber\\
&=& \frac{a_{j+n}}{a_j}\times \frac{\displaystyle\prod_{i=1}^{n-1} (i+j+\gamma-1)}{\displaystyle\prod_{i=1}^{n-1} (i+j+\gamma+n-1)} \nonumber\\
&=& \frac{\displaystyle\prod_{i=1}^{n} \bigl[\alpha(i+j+\gamma-1)+\beta\bigr]}{\displaystyle\prod_{i=1}^{n} (i+j+\gamma-1)}\times \frac{\displaystyle\prod_{i=1}^{n-1} (i+j+\gamma-1)}{\displaystyle\prod_{i=1}^{n-1} (i+j+\gamma+n-1)} \nonumber\\
&=& \prod_{i=1}^{n} \frac{\alpha(i+j+\gamma-1)+\beta}{i+j+\gamma+n-2}.
\end{eqnarray}

Joining the last result with Remark \ref{Remark_d_n_k_null}, we conclude that the relation
\begin{equation} \label{Recurrence_d_n_k}
d_n^{(j+1)} = d_n^{(j)} \prod_{i=1}^{n} \frac{\alpha(i+j+\gamma-1)+\beta}{i+j+\gamma+n-2} 
\end{equation}
is always valid and consequently, we obtain the formula
\begin{equation} \label{Cute_formula_d_n_k}
d_n^{(k)} = d_n^{(0)}\ \prod_{j=0}^{k-1} \prod_{i=1}^{n} \frac{\alpha(i+j+\gamma-1)+\beta}{i+j+\gamma+n-2}.
\end{equation}

Now let us focus on $d_n^{(0)}$. Following Equation \eqref{Ugly_formula_d_n_k}, we can write
\begin{equation} \label{Ugly_formula_d_n_0}
d_n^{(0)} = a_0 a_{1} \ldots a_{n-1}\ \frac{\displaystyle\prod_{i=1}^{n-1} i! \: \bigl[\alpha(i-1) - \beta \bigr]^{n-i}}{\displaystyle\prod_{i=1}^{n-1} \bigl[i+\gamma-1\bigr]^i \: \bigl[2n-i+\gamma-2\bigr]^i}.
\end{equation}

It is easy to see that similarly to Remark \ref{Remark_d_n_k_null}, $d_n^{(0)} = 0$ implies $d_{n+1}^{(0)} = 0$. Again, if we have $d_j^{(0)} \neq 0$ for some $j \geq 1$, we compute the ratio
\begin{eqnarray}
\frac{d_{j+1}^{(0)}}{d_j^{(0)}} &=& a_j\: j!\: \prod_{i=1}^{j} \bigl[\alpha(i-1)-\beta\bigr]\times \frac{\displaystyle\prod_{i=1}^{j-1} \bigl[i+\gamma-1\bigr]^i \: \bigl[2j-i+\gamma-2\bigr]^i}{\displaystyle\prod_{i=1}^{j} \bigl[i+\gamma-1\bigr]^i \: \bigl[2j-i+\gamma\bigr]^i} \nonumber\\
&=& a_j\: j! \ \times \frac{\displaystyle\prod_{i=1}^{j} \bigl[\alpha(i-1)-\beta\bigr]}{\displaystyle\prod_{i=1}^{j} \bigl[i+j+\gamma-2\bigr] \bigl[i+j+\gamma-1\bigr]} \nonumber\\
&=& j!\ \prod_{i=1}^{j} \frac{\alpha(i+\gamma-1) + \beta}{i+\gamma-1} \times \frac{\displaystyle\prod_{i=1}^{j} \bigl[\alpha(i-1)-\beta\bigr]}{\displaystyle\prod_{i=1}^{j} \bigl[i+j+\gamma-2\bigr] \bigl[i+j+\gamma-1\bigr]} \nonumber\\
&=& \prod_{i=1}^{j} \frac{i\bigl[\alpha(i+\gamma-1) + \beta\bigr]\bigl[\alpha(i-1)-\beta\bigr]}{\bigl[ i+ \gamma -1 \bigr]\bigl[i+j+\gamma-2\bigr]  \bigl[i+j+\gamma-1\bigr]}.
\end{eqnarray}

Hence, we can conclude that the relation
\begin{equation}
d_{j+1}^{(0)} = d_j^{(0)}\:\prod_{i=1}^{j} \frac{i\bigl[\alpha(i+\gamma-1) + \beta\bigr]\bigl[\alpha(i-1)-\beta\bigr]}{\bigl[i+\gamma-1 \bigr]\bigl[i+j+\gamma-2\bigr]  \bigl[i+j+\gamma-1\bigr]}
\end{equation}
is always valid and therefore, from the obvious $d_1^{(0)} = a_0 = 1$ we obtain
\begin{equation} \label{Cute_formula_d_n_0}
d_n^{(0)} = \prod_{1 \leq i \leq j \leq n-1} \frac{i\bigl[\alpha(i+\gamma-1) + \beta\bigr]\bigl[\alpha(i-1)-\beta\bigr]}{\bigl[i+\gamma-1 \bigr]\bigl[i+j+\gamma-2\bigr]  \bigl[i+j+\gamma-1\bigr]}.
\end{equation}

We summarize the obtained results in the following theorem.

\begin{Theorem} \label{Theorem_Hankel_transform_a_n}
The generalized Hankel transform of the sequence $(a_n)_{n=0}^{\infty}$ with $a_0 = 1$ and satisfying \eqref{Rel_def_a_n} is given by
\begin{multline} \label{Principal_Identity}
d_n^{(k)} = \prod_{1 \leq i \leq j \leq n-1} \frac{i\bigl[\alpha(i+\gamma-1) + \beta\bigr]\bigl[\alpha(i-1)-\beta\bigr]}{\bigl[i+\gamma-1 \bigr]\bigl[i+j+\gamma-2\bigr]  \bigl[i+j+\gamma-1\bigr]} \\
\times \prod_{j=0}^{k-1} \prod_{i=1}^{n} \frac{\alpha(i+j+\gamma-1)+\beta}{i+j+\gamma+n-2}.
\end{multline}
Moreover, if we allow $a_0$ to be arbitrary, then the expression for $d_n^{(k)}$ gets multiplied by $a_0^n$.
\end{Theorem}

\section{Applications to Particular Sequences}
\label{SecnAppl}

We now consider several applications of Theorem \ref{Theorem_Hankel_transform_a_n} to particular important sequences. In some cases, we show how Formula \eqref{Principal_Identity} applies to derive some well known closed product form evaluations, which have been obtained before using several independent methods. Bearing this goal in mind, we do not prove most product identities; the avid reader will easily be able to supply the proofs.

We recall that in this paper, empty products are always considered to be $1$.

\begin{Example}
Let $a_n = (n+\kappa)^{-1}$, whose associated matrices $H_n^{(k)}$ are (generalized) Hilbert matrices. The sequence $(a_n)$ satisfies the hypotheses of Theorem \ref{Theorem_Hankel_transform_a_n} for $\alpha=1$, $\beta=-1$, $\gamma=1+\kappa$ and thus, its Hankel transform is given by
\begin{multline}
d_n^{(k)} = \frac{1}{\kappa^n}\ \prod_{1 \leq i \leq j \leq n-1} \frac{i^2 \bigl[i-1+\kappa  \bigr]}{\bigl[i+\kappa \bigr]\bigl[i+j+\kappa-1\bigr]  \bigl[i+j+\kappa\bigr]}\\ 
\times \prod_{j=0}^{k-1} \prod_{i=1}^{n} \frac{i+j+\kappa-1}{i+j+n+\kappa-1}.
\end{multline}

After some elementary transformations, we obtain the identities
\begin{eqnarray*}
\prod_{1 \leq i \leq j \leq n-1} \frac{i^2 \bigl[i-1+\kappa  \bigr]}{\bigl[i+\kappa \bigr]\bigl[i+j+\kappa-1\bigr]  \bigl[i+j+\kappa\bigr]} &=& \kappa^n\ \frac{\displaystyle\prod_{i=1}^{n-1} (i!)^2}{\displaystyle\prod_{i,j=1}^n (i+j+\kappa-2)}, \\
\prod_{j=0}^{k-1} \prod_{i=1}^{n} \frac{i+j+\kappa-1}{i+j+n+\kappa-1} &=& \prod_{i,j=1}^n \frac{i+j+\kappa-2}{i+j+k+\kappa-2};
\end{eqnarray*}
which allow us to write $d_n^{(k)}$ in the more familiar form
\begin{equation}
d_n^{(k)} = \frac{\displaystyle\prod_{i=1}^{n-1} (i!)^2}{\displaystyle\prod_{i,j=1}^n (i+j+k+\kappa-2)}.
\end{equation}

This is a very well known formula (especially for $\kappa=1$) and can be proved by several methods (for some historical remarks, see Muir \cite[vol. III, pp. 311]{Muir}). For example, it can be easily derived from \emph{Cauchy's double alternant} (see Krattenthaler \cite[Eq. 2.7]{Kratt1}). An extensive literature exists on Hilbert matrices and their generalizations. For an interesting compilation of results about Hilbert matrices, please refer to Choi \cite{Choi}. For a study from the viewpoint of orthogonal polynomials, see Berg \cite{Berg}.
\end{Example}


\begin{Example}
 Let $a_n = 2(n^2+3n+2)^{-1}$ be the sequence of the reciprocals of triangular numbers. It satisfies \eqref{Rel_def_a_n}  for $\alpha=1$, $\beta=-2$, $\gamma=3$. Thus, its Hankel transform is given by
\begin{equation}
d_n^{(k)} = \prod_{1 \leq i \leq j \leq n-1} \frac{i^2\bigl[i+1\bigr]}{\bigl[i+2 \bigr]\bigl[i+j+1\bigr]  \bigl[i+j+2\bigr]} \times \prod_{j=0}^{k-1} \prod_{i=1}^{n} \frac{i+j}{i+j+n+1}.
\end{equation}

By considering the identities
\begin{eqnarray*}
\prod_{1 \leq i \leq j \leq n-1} \frac{i^2\bigl[i+1\bigr]}{\bigl[i+2 \bigr]\bigl[i+j+1\bigr]  \bigl[i+j+2\bigr]} &=& 2^n\ \frac{\displaystyle\prod_{i=1}^{n-1} (i!)^2}{\displaystyle\prod_{i,j=1}^n (i+j)}, \\
\prod_{j=0}^{k-1} \prod_{i=1}^{n} \frac{i+j}{i+j+n+1} &=& \binom{n+k}{n}^{-1} \prod_{i,j=1}^n \frac{i+j}{i+j+k};  
\end{eqnarray*}
we obtain the simpler formula
\begin{equation}
d_n^{(k)} = 2^n \binom{n+k}{n}^{-1} \times \frac{\displaystyle\prod_{i=1}^{n-1} (i!)^2}{\displaystyle\prod_{i,j=1}^n (i+j+k)}.
\end{equation}
\end{Example}

\begin{Remark}
Consider the apparently similar sequence $a_n = (n+1)^{-2}$. Clearly, it does not satisfy \eqref{Rel_def_a_n} for any $\alpha, \beta, \gamma$. Actually, its Hankel transform is unlikely to have a closed product form evaluation. Indeed, note the factorizations
\begin{eqnarray}
d_3^{(0)} &=& \frac{647}{2^{8} \cdot 3^6 \cdot 5^2}, \\
d_5^{(0)} &=& \frac{179 \cdot 179357}{2^{20} \cdot 3^{6} \cdot 5^{10} \cdot 7^{5}}, \\
d_7^{(0)} &=& \frac{23\cdot 1280587616051046200369}{2^{36} \cdot 3^{22} \cdot 5^{10}\cdot 7^{14} \cdot 11^{6} \cdot 13^{2}}.
\end{eqnarray}
The amazingly large primes in the numerators suggest the claim.
\end{Remark}


\begin{Example}
Consider the sequence $a_n = (n!)^{-1}$. It satisfies the hypotheses of Theorem \ref{Theorem_Hankel_transform_a_n} for $\alpha=0$, $\beta=1$, $\gamma=1$ and therefore, we have
\begin{equation}
d_n^{(k)} = (-1)^{\binom{n}{2}} \prod_{1 \leq i \leq j \leq n-1} \frac{1}{\bigl[i+j-1\bigr]  \bigl[i+j\bigr]} \times \prod_{j=0}^{k-1} \prod_{i=1}^{n} \frac{1}{i+j+n-1}.
\end{equation}

Alternatively, the following identities hold:
\begin{eqnarray*}
\prod_{1 \leq i \leq j \leq n-1} \frac{1}{\bigl[i+j-1\bigr]  \bigl[i+j\bigr]} &=& \prod_{i=0}^{n-1} \frac{i!}{(i+n-1)!}, \\
\prod_{j=0}^{k-1} \prod_{i=1}^{n} \frac{1}{i+j+n-1} &=& \prod_{i=0}^{n-1} \frac{(i+n-1)!}{(i+k+n-1)!}. 
\end{eqnarray*}

Hence, we obtain
\begin{equation}
d_n^{(k)} = (-1)^{\binom{n}{2}} \prod_{i=0}^{n-1} \frac{i!}{(i+k+n-1)!}, 
\end{equation}
which recovers the formula from Bacher \cite[Thm. 1.3]{Bacher}. 
\end{Example}


\begin{Example}
Let $a_n = (n+1)^{-1} \binom{2n}{n}$ be the sequence of Catalan numbers. It satisfies the hypotheses of Theorem \ref{Theorem_Hankel_transform_a_n} for $\alpha=4$, $\beta=-6$ and $\gamma=2$. Therefore, its Hankel transform is given by
\begin{equation}
d_n^{(k)} = \prod_{1 \leq i \leq j \leq n-1} \frac{i\bigl[4i-2\bigr]\bigl[4i+2\bigr]}{\bigl[i+1 \bigr]\bigl[i+j\bigr]  \bigl[i+j+1\bigr]} \times \prod_{j=0}^{k-1} \prod_{i=1}^{n} \frac{4(i+j)-2}{i+j+n}.                                                                                                                                                                                                                                                                                                                                                                                                                                  
\end{equation}

The left product reduces to $1$ and the right product can be rewritten as
\begin{equation*}
\prod_{1\leq i \leq j \leq k-1} \frac{i+j+2n}{i+j}
\end{equation*}
for all $k \geq 0$, which is obvious from the identity
\begin{equation*}
\prod_{i=1}^{n} \frac{4(i+j)-2}{i+j+n} = \prod_{i=1}^{j} \frac{i+j+2n}{i+j},\qquad j \geq 0.
\end{equation*}

Accordingly, we obtain the well known formula
\begin{equation} \label{Formula_Catalan_Numbers}
d_n^{(k)} = \prod_{1\leq i \leq j \leq k-1} \frac{i+j+2n}{i+j},
\end{equation}
which was primarily found by Desainte-Catherine and Viennot \cite{DCV}, who also gave a combinatorial interpretation for this transform (see also Gessel and Viennot \cite{GV} for further generalizations). It is also proved in Tamm \cite{Tamm}, by means of the Dodgson's condensation method (see Krattenthaler \cite[Prop. 10]{Kratt1} for details) and discussed in Krattenthaler \cite[Thm. 33]{Kratt2}, with some additional remarks. The cases $k=0, 1$ and $2$ are also studied in Aigner \cite{Aigner}, where the author describes a beautiful generalization of Catalan numbers (called Catalan -- like numbers) inspired by the property $d_n^{(0)} = d_n^{(1)}=1$.
\end{Example}


\begin{Example}
Let $a_n = \binom{2n}{n}$ be the sequence of even central binomial coefficients. It is immediate to see that $(a_n)$ satisfies the hypotheses of Theorem \ref{Theorem_Hankel_transform_a_n} for $\alpha=4$, $\beta=-2$ and $\gamma=1$. Hence, its Hankel transform is given by
\begin{equation}
d_n^{(k)} = \prod_{1 \leq i \leq j \leq n-1} \frac{\bigl[4i - 2 \bigr]^2}{\bigl[i+j-1\bigr]  \bigl[i+j\bigr]} \times \prod_{j=0}^{k-1} \prod_{i=1}^{n} \frac{4(i+j)-2}{i+j+n-1}.
\end{equation}

The left product is readily seen to be $2^{n-1}$. As for the right product, it can be rewritten as
\begin{equation*}
2 \times \prod_{1 \leq i \leq j \leq k-1} \frac{i+j-1+2n }{ i+j-1}
\end{equation*}
for all $k \geq 1$. This can be deduced directly from the identities
\begin{eqnarray*}
\prod_{i=1}^{n} \frac{4(i+j)-2}{i+j+n-1} &=& \prod_{i=1}^{j} \frac{i+j-1+2n}{i+j-1},\qquad j\geq 1,\\
\prod_{i=1}^{n} \frac{4i-2}{i+n-1} &=& 2.
\end{eqnarray*}

Thus we are able to obtain the formula
\begin{equation}
d_n^{(k)} = \begin{cases}
2^{n-1}, & \textrm{if $k=0$;} \\
2^n\times \displaystyle\prod_{1 \leq i \leq j \leq k-1} \frac{i+j-1+2n }{i+j-1},  & \textrm{if $k \geq 1$.}
\end{cases}
\end{equation}

Taking into account that $\binom{2m}{m} = 2\binom{2m-1}{m}$, we can recover the formula for the Hankel transform of odd central binomial coefficients from Tamm \cite[Eq. 1.5]{Tamm}. It is worth mentioning that, given the identity $\binom{2n}{n} = (-1)^n\: 2^{2n}\: \binom{-1/2}{n}$, the Hankel transform could also be computed directly using \cite[Thm. 26, Eq. 3.12]{Kratt1}. For interesting connections of this Hankel transform with combinatorics and algebra, see Aigner \cite{Aigner} and Garcia Armas and Sethuraman \cite{Garcia_Seth}. 
\end{Example}

It is worth mentioning that several other interesting sequences satisfy \eqref{Rel_def_a_n} and their Hankel transforms can therefore be evaluated using Theorem \ref{Theorem_Hankel_transform_a_n}, e.g., the binomial sequences $a_n = \binom{\lambda}{n}$, where $\lambda \in \mathbb{C}$, and 
$b_n = \binom{n + \lambda}{m}$, where $m \in \mathbb{Z}_{\geq 0}$ and $\lambda \in \mathbb{Z}$ with $\lambda \geq m$, or $\lambda \in \mathbb{C} \backslash \mathbb{Z}$.

\subsection{Hankel Transforms of Reciprocals}
We finish this section by noting the following beautiful property: if $(a_n)$ is a non-zero sequence satisfying Equation \eqref{Rel_def_a_n} with $\alpha \neq 0$, then the reciprocal sequence $(a_n^{-1})$ satisfies the relation
\begin{equation} \label{Def_rel_reciprocal_seq}
a_{n+1}^{-1} = \frac{n+\gamma}{\alpha n +\alpha \gamma + \beta}\; a_n^{-1} = \biggl(\frac{1}{\alpha} - \frac{\frac{\beta}{\alpha^2}}{n+\gamma+\frac{\beta}{\alpha}}\biggr)\ a_n^{-1}.
\end{equation}

\begin{Corollary} \label{Cor_Reciprocal_Sequence}
Let $(a_n)_{n=0}^{\infty}$ be a non-zero sequence with $a_0=1$ and satisfying \eqref{Rel_def_a_n} for some $\alpha \neq 0$. Then, the generalized Hankel transform $d_n^{(k)}$ of the reciprocal sequence $(a_n^{-1})_{n=0}^{\infty}$ is given by
\begin{multline} \label{Transform_Reciprocal_Seq}
\prod_{1 \leq i \leq j \leq n-1} \frac{i\bigl[i+\gamma-1\bigr]\bigl[\alpha(i-1)+\beta\bigr]}{\bigl[\alpha (i+\gamma-1)+\beta\bigr]\bigl[\alpha(i+j+\gamma-2)+\beta\bigr]  \bigl[\alpha(i+j+\gamma-1)+\beta\bigr]} \\
\times \prod_{j=0}^{k-1} \prod_{i=1}^{n} \frac{i+j+\gamma-1}{\alpha(i+j+\gamma+n-2)+\beta}.
\end{multline}
\end{Corollary}

\begin{Remark}
The formula remains valid in the case $\alpha = 0$; an easy way to see this is by making $\alpha \to 0$ in \eqref{Transform_Reciprocal_Seq}.
\end{Remark}

As an immediate consequence of Corollary \ref{Cor_Reciprocal_Sequence}, we evaluate some generalized Hankel transforms:
\begin{itemize}
\item Let $a_n = (n+1) \binom{2n}{n}^{-1}$ be the sequence of the reciprocals of Catalan numbers. Then, its Hankel transform is given by
\begin{multline}
d_n^{(k)} = \frac{1}{2^{n(n+k-1)}}\ \prod_{1 \leq i \leq j \leq n-1} \frac{i\bigl[i+1\bigr]\bigl[2i-5\bigr]}{\bigl[2i-1\bigr]\bigl[2(i+j)-3\bigr]  \bigl[2(i+j)-1\bigr]} \\
\times \prod_{j=0}^{k-1} \prod_{i=1}^{n} \frac{i+j+1}{2(i+j+n)-3}.
\end{multline}

\item Let $a_n = \binom{2n}{n}^{-1}$ be the sequence of the reciprocals of even central binomial coefficients. Then, its Hankel transform is given by
\begin{multline}
d_n^{(k)} = \frac{1}{2^{n(n+k-1)}}\ \prod_{1 \leq i \leq j \leq n-1} \frac{i^2\bigl[2i-3\bigr]}{\bigl[2i-1\bigr]\bigl[2(i+j)-3\bigr]  \bigl[2(i+j)-1\bigr]} \\
\times \prod_{j=0}^{k-1} \prod_{i=1}^{n} \frac{i+j}{2(i+j+n)-3}.
\end{multline}
\end{itemize}

The form of the above determinants, together with extensive computational evidence collected by the author, suggest the following conjecture.

\begin{Conjecture}
Let $(a_n)_{n=0}^{\infty}$ be the sequence of the reciprocals of Catalan numbers, or the sequence of the reciprocals of even central binomial coefficients. Then, the Hankel matrices $H_n^{(k)}$ associated to $(a_n)$ have inverses whose entries are all integers.
\end{Conjecture}

\section*{Acknowledgments} The author would like to thank Prof. B. A. Sethuraman for introducing him to the subject and for so many fruitful comments and discussions. The author also thanks C. Krattenthaler for pointing out to him how to derive the main formula in this paper from the results in \cite{Kratt1}.


\begin{thebibliography}{99}
\bibitem{Aigner} M. Aigner, Catalan--like numbers and determinants,
\textit{J. Combin. Theory Ser. A} \textbf{87}  (1999), 33--51.
\bibitem{Bacher} R. Bacher, Determinants of matrices related to the Pascal triangle, \emph{J. Th\' eor. Nombres Bordeaux} \textbf{14} (2002), 19--41.
\bibitem{Berg} C. Berg, Fibonacci numbers and orthogonal polynomials, arXiv:math/0609283v2, preprint.
\bibitem{Chamb_French} M. Chamberland and C. French, Generalized Catalan numbers and generalized Hankel transformations, \textit{J. Integer Seq.} \textbf{10} (2007), Article 07.1.1.
\bibitem{Choi} M. Choi, Tricks or treats with the Hilbert matrix, \emph{Amer. Math. Monthly} \textbf{80} (1983), 301--312.
\bibitem{CRI} A. Cvetkovi\'c, P. Rajkovi\'c and M. Ivkovi\'c, Catalan numbers, the Hankel transform, and Fibonacci numbers, \emph{J. Integer Seq.} \textbf{5} (2002), Article 02.1.3.
\bibitem{DCV} M. Desainte-Catherine and X. G. Viennot, Enumeration of certain Young tableaux with bounded height, \emph{Combinatoire \' Enum\' erative (Montreal 1985)}, Lect. Notes in Math. \textbf{1234} (1986), 58--67 
\bibitem{ERR_Riemann} \"{O}. Egecioglu, T. Redmond and C. Ryavec, From a polynomial Riemann hypothesis to alternating sign matrices, \emph{Electron.
J. Combin.} \textbf{8} R36 (2001).
\bibitem{ERR} \"{O}. Egecioglu, T. Redmond and C. Ryavec,
Almost product evaluation of Hankel determinants, \emph{Electron.
J. Combin.} \textbf{15}(1)\ R6 (2008).
\bibitem{French} C. French, Transformations preserving the Hankel transform, \emph{J. Integer Seq.} \textbf(10) (2007), Article 07.7.3.
\bibitem{Garcia_Seth} M. Garcia Armas and B. A. Sethuraman, A note on the Hankel transform of the central binomial coefficients, \emph{J. Integer Seq.} \textbf{11} (2008), Article 08.5.8.
\bibitem{GV} I. M. Gessel and X. G. Viennot, Determinants, paths, and plane partitions, preprint (1989).
\bibitem{Kratt1} C. Krattenthaler, Advanced determinant calculus, \emph{S\' em. Lothar. Combin.} \textbf{42} (1999), Article B42q.
\bibitem{Kratt2} C. Krattenthaler, Advanced determinant calculus: a complement, \emph{Linear Algebra Appl.} \textbf{411} (2005), 68--166.
\bibitem{Lay} J. W. Layman, The Hankel transform and some of
its properties, \emph{J. Integer Seq.} \textbf{4} (2001), Article 01.1.5.
\bibitem{Muir} T. Muir, \emph{The Theory of Determinants in the Historical Order of Development}, 4 vols., MacMillan, 1906--1923.
\bibitem{SpSt} M. Z. Spivey and L. L. Steil, The $k$-binomial
transform and the Hankel transform, \emph{J. Integer Seq.} \textbf{9} (2006), Article 06.1.1.
\bibitem{Tamm} U. Tamm, Some aspects of Hankel matrices in coding theory
and combinatorics, \emph{Electron. J. Combin.}
\textbf{8}(1)\ A1 (2001).
\bibitem{Vein_Dale} R. Vein and A. Dale, \emph{Determinants and Their Applications in Mathematical Physics}, Springer, 1991.
\end{thebibliography}
\end{document}